\documentclass[12 pt]{article}
\usepackage[T1]{fontenc}
\usepackage[utf8x]{inputenc}
\usepackage{lmodern}
\usepackage{amsmath,amsthm,amsfonts,amssymb}    %Package (presque) indispensable pour taper corectement des maths
\usepackage{graphicx}
\usepackage{enumitem}
\usepackage[margin=25mm]{geometry}  %Pour avoir des marges plus petites
\usepackage{hyperref}   %Pour que les equations, enoncés référencés ou sections dans la table des matières soient cliquables
\usepackage{stmaryrd}
\usepackage{picture}
\usepackage{mathtools}
\usepackage{tikz}   %Pour faire des figures
\usepackage[french,intoc]{nomencl}    %Pour pouvoir générer une liste de symboles
\makenomenclature   %Pour générer une liste de symboles en fin de mémoire
\usepackage{makeidx}
\makeindex
\usepackage{bbold}  %Pour faire de belles fonctions indicatrices
\usepackage{pgfplots}   %pour les figures Tikz générées par Geogebra
\pgfplotsset{compat=1.15}   %pour les figures Tikz générées par Geogebra
\usepackage{mathrsfs}   %pour les figures Tikz générées par Geogebra
\usetikzlibrary{arrows}   %pour les figures Tikz générées par Geogebra
\usepackage{subfigure}
\usepackage{xcolor}
\usepackage{eurosym}
\usepackage{fourier}
\usepackage{tkz-tab}    %Pour les tableaux de signes et variations
\usepackage{diagbox}    %Pour couper une case d'un tableau par une diagonale
\usepackage{pifont}     %Pour utiliser des icônes
\allowdisplaybreaks     %Pour que les align puissent être sur deux pages consécutives et donc que ça ne laisse pas de gros blancs
\usepackage{authblk}

\title{Numerical range of finite perturbations of a 2-periodic weighted shift}
\author{Simon Marionnet \thanks{simon.marionnet@univ-lille.fr}}
\date{2026}
\affil{Univ. Lille, CNRS, UMR 8524 - Laboratoire Paul Painlevé, F-59000 Lille, France}

\newcommand{\N}{\mathbb{N}}

\newcommand{\R}{\mathbb{R}}

\newcommand{\C}{\mathbb{C}}

\newcommand{\D}{\mathbb{D}}

\DeclareMathOperator{\Tr}{Tr}
\DeclareMathOperator{\RE}{Re}

\DeclareMathOperator{\com}{com}

\theoremstyle{definition}
\newtheorem{defi}{Definition}[section]  %On numérote les résultats selon la section.
\newtheorem{prop}[defi]{Proposition}    %On fait en sorte que chaque numéro ne renvoie qu'a un seul résultat. Par exemple il ne peut pas y avoir une définition 1.3.1 et un théorème 1.3.1
\newtheorem{thm}[defi]{Theorem}
\newtheorem{lem}[defi]{Lemma}
\newtheorem{coro}[defi]{Corollary}

\theoremstyle{remark}
\newtheorem{rem}[defi]{Remark}
\newtheorem{ex}[defi]{Example}

\newcommand{\subjclass}[1]{%
  \par\noindent\textbf{2020 Mathematics Subject Classification:} #1\par
}

\newcommand{\keywords}[1]{%
  \par\noindent\textbf{Keywords:} #1\par
}

\begin{document}

\maketitle

\begin{abstract}
    In this note, we consider a weighted shift operator $S_\alpha$ of weight sequence\\
    $\alpha=(a_1,\dots,a_n,a,b,a,b,\dots)$. We characterise the equality of the closure of the numerical range of a $S_\alpha$ and its essential numerical range. We end up with some applications and comments.
\end{abstract}

\keywords{Numerical range, essential numerical range, weighted shift operator}

\subjclass{47A12, 47B37}

\section{Introduction}

For a bounded operator $T$ on a complex Hilbert space $H$, the \textit{numerical range} of $T$ is
\begin{equation*}
    W(T):=\left\{\langle Tx,x\rangle: x\in H \text{ and }\lVert x\rVert=1\right\}\subseteq \C.
\end{equation*}
The \textit{numerical radius} of $T$ is $w(T):=\sup\left\{\lvert z\rvert: z\in W(T)\right\}$. Cauchy-Schwarz inequality leads to $w(T)\leq \lVert T\rVert$. From the definition we have that $W(\lambda T)=\lambda W(T)$ for all $\lambda\in\C$.

It is well known that for every operator $T$, the set $W(T)$ is convex and the following inclusions hold $\sigma_p(T)\subseteq W(T)$ and $\sigma(T)\subseteq\overline{W(T)}$ where $\sigma_p(T)$ and $\sigma(T)$ are respectively the point spectrum (\textit{i.e.} the set of eigenvalues) and the spectrum of $T$. Furthermore, $W(T)=W(U^*TU)$ for every unitary operator $U$. $T$ and $U^*TU$ are said to be \textit{unitary equivalent}. One may find all useful material about the numerical range and the numerical radius of an operator in \cite{gau2021numerical}.

Throughout the paper, the \textit{real part} of $T$ will be denoted by $\RE(T):=\frac{T+T^*}{2}$. We have $W(\RE(T))=\RE(W(T))$.

A self-adjoint operator $T\in\mathcal{L}(H)$ is called \textit{positive} if $\langle Tx,x\rangle\geq 0;\forall x\in H$. We write $A\leq B$ to express the fact that $B-A\geq0$.

The \textit{essential numerical range} of an operator $T$ is
\begin{equation*}
    W_e(T):=\bigcap\left\{\overline{W(T+K)}: K \text{ is a compact operator}\right\}.
\end{equation*}
The \textit{essential numerical radius} of $T$ is $w_e(T):=\sup\left\{\lvert z\rvert: z\in W_e(T)\right\}$. It is clear that $W_e(T)$ is closed and $W_e(T)\subseteq \overline{W(T)}$. A result of Lancaster states that if $W(T)$ has no extreme point, $W_e(T)=\overline{W(T)}$ \cite[Corollary 2]{lancaster1975boundary}. In particular, this holds when $W(T)$ is open.

Let $K$ be a closed subspace of $H$ and $A,B$ bounded operators respectively on $H$ and $K$. $B$ is said to be a \textit{compression} of $A$, or $A$ is a \textit{dilation} of $B$, if $B=\pi_K A\pi_K$ where $\pi_K$ is the orthogonal projection from $H$ onto $K$. If $B$ is a compression of $A$, we have $W(B)\subseteq W(A)$.

For a bounded sequence $\alpha=(\alpha_n)_{n\geq0}\in \ell^\infty$, the unilateral \textit{weighted shift} of weight sequence $\alpha$ denoted by $S_\alpha$ maps a sequence $(x_0,x_1,\dots)\in \ell^2$ to $(0,\alpha_0x_0,\alpha_1x_1,\dots)$. That defines a bounded operator on the Hilbert space $\ell^2$ and $\lVert S_\alpha\rVert=\lVert \alpha\rVert_\infty$. For a sequence $\alpha=(\alpha_n)_n\in\ell^\infty$, let us denote by $|\alpha|:=(|\alpha_n|)_n$. It is well known that $S_\alpha$ is unitary equivalent to $S_{|\alpha|}$. Thus, for all $\alpha\in \ell^\infty$ and $\theta\in [0,2\pi]$, $e^{i\theta}S_\alpha$ is unitary equivalent to $S_\alpha$. Hence, the numerical range of a weighted shift is a disk centred at 0 by convexity. Then, to completely describe the numerical range of a weighted shift, we just need to find the radius of the disk and to determine if it is open or closed. For example, the classical shift $S$ (the weighted shift with weight sequence $\alpha=(1,1,\dots)$) has the open unit disk $\D$ as numerical range.

Since $S_\alpha$ is unitary equivalent to $S_{|\alpha|}$ we will only consider non-negative weight sequences in the following. If there is a $k$ such that $\alpha_k=0$, then $S_\alpha=S_1\oplus S_2$ where $S_1$ and $S_2$ are (maybe finite) weighted shifts. Since the numerical range of a direct sum is the convex hull of the union of the numerical ranges, we can restrict ourselves to positive weight sequences.

In the following, we will denote by $D(0,r)$ and $\overline{D}(0,r)$ respectively the open and closed disk of radius $r\geq 0$ centred at 0.

The numerical range of weighted shifts with periodic weight sequence was first studied by Ridge in \cite{ridge1976numerical}. He proved that a weighted shift with periodic weight sequence \\
$(a_1,\dots,a_n,a_1,\dots,a_n,a_1,\dots)$ has the same numerical radius as the weighted cyclic matrix \\$\begin{pmatrix}
    0&&&a_n\\
    a_1&0&&\\
    &\ddots&\ddots&\\
    &&a_{n-1}&0
\end{pmatrix}$. Thus, for $n=2$, we know that $w(S_\alpha)=(a+b)/2$ if $\alpha=(a,b,a,b,a\dots)$. Nevertheless, we do not know the numerical radius of periodic shifts with period bigger than 3. Indeed, in dimension bigger than 3, the numerical radius of weighted cyclic matrices is unknown. This is why we will focus on 2-periodic weighted shifts in the sequel.

Ridge conjectured that the numerical range of periodic weighted shifts are open. That was settled in positive by Stout in \cite{stout1983numerical} (see also \cite[Theorem 4]{tam1997conjecture} and \cite[Proposition 4.2]{WangWu2011NRofWS}). Thus, we can claim that for every periodic sequence $\alpha$, $W_e(S_\alpha)=\overline{W(S_\alpha)}$ by Lancaster's Theorem cited above. Fix $n$ an even number, $a_1,\dots,a_n>0$ and $a,b>0$ and define $\alpha:=(a_1,\dots,a_n,a,b,a,b,\dots)$. $S_\alpha$ is a finite (thus compact) perturbation of the 2-periodic shift with weight sequence $(a,b,a,b,\dots)$.

According to the preceding, the situation is as follows:
\begin{enumerate}
    \item $W(S_\alpha)$ is a disk centred at 0 with radius $(a+b)/2=w_e(S_\alpha)$;
    \item or $W(S_\alpha)$ is a closed disk centred at 0 with radius $r>(a+b)/2$.
\end{enumerate}
Theorem 4.3 of \cite{WangWu2011NRofWS} ensures that in the first case, $W(S_\alpha)$ is open. We will give another proof of this result at the end of the article in our context of 2-periodic shifts. Our goal in this article is to give a characterisation of when the first case appears. That is our Theorem \ref{central2}. We will follow the philosophy of the article of Wang and Wu \cite{WangWu2011NRofWS}, and will not try to express the numerical radius in the second case.

This question has already been studied in the past. In 1967, Berger and Stampfli treated the case $n=1$ and $a=b=1$ in \cite{berger1967mapping}. If $t>0$ and $\alpha=(t,1,1,\dots)$, $W(S_\alpha)=\D$ if $t\leq \sqrt{2}$ and equals $\overline{D}\left(0,\frac{t^2}{2\sqrt{t^2-1}}\right)$ if $t>\sqrt{2}$. The proof uses Hardy space techniques. The fact that $W(S_\alpha)=\D$ is shown to be equivalent for the polynomial $(t^2-1)z^2-1$ to have no root in $(0,1)$.

This approach has been followed in \cite[Corollary 4.8]{WangWu2011NRofWS} (see also \cite{chien2013numerical}) and \cite{vandanjav2014numerical} respectively for the weight sequence $\alpha=(1,h,1,1\dots)$ and $\alpha=(a_1,a_2,1,1,\dots)$. We have $w(S_\alpha)=1$ if and only if $h\leq \sqrt{6}/2$ and $4-a_1^2-2a_2^2\geq0$ respectively. Otherwise, the numerical radius of $S_\alpha$ is provided.

The case of an arbitrary finite perturbation $\alpha=(a_1,\dots,a_n,1,1,\dots)$ was first considered in \cite{undrakh2015numerical}. The authors used the same method as above and exhibit a polynomial $G_\alpha$ such that $w(S_\alpha)>1$ if and only if $G_\alpha$ has a real root in $(0,1)$ (Theorem 3.3 of \cite{undrakh2015numerical}). The numerical radius is no longer explicit there, but is expressed in terms of the smallest root of $G_\alpha$ in $(0,1)$.

These results have been extended to the case of a finite perturbation of a 2-periodic weighted shift. The case of $\alpha=(a_1,a_2,a,b,a,b,\dots)$ has been studied in \cite{chakraborty2022numerical}.

The case $\alpha=(a_1,\dots,a_n,a,b,a,b,\dots)$ was the object of \cite{chakraborty2022numericalradii}. Remark 2.7 of the latter ensures that $w(S_\alpha)>w_e(S_\alpha)=(a+b)/2$ if and only if the polynomial
\begin{equation*}
    G_{\alpha,e}(z)=abQ_{n-1,\alpha}(z)-a_n^2(a^2z^2+ab)z^2\frac{Q_{n-2,\alpha}(z)}{\sqrt{abz^4+(a^2+b^2)z^2+ab}}
\end{equation*}
has a real root between 0 and 1 if $n$ is even. If $n$ is odd, thanks to Theorem 2.6 of the same article, we have the same criteria with the polynomial
\begin{equation*}
    G_{\alpha,o}(z)=(a^2z^2+ab)\frac{Q_{n-1,\alpha}(z)}{\sqrt{abz^4+(a^2+b^2)z^2+ab}}-a_n^2z^2Q_{n-2,\alpha}(z).
\end{equation*}
The polynomial $Q_{k,\alpha}(z)$ is defined by
\begin{equation*}
    Q_{k,\alpha}(z)=\det\left(\sqrt{abz^4+(a^2+b^2)z^2+ab}I_{k+1}+2z\RE\left(A_{\alpha_{|k}}\right)\right),
\end{equation*}
where $\alpha_{|k}$ is the sequence built of the $k$ first terms of the sequence $\alpha$. Initially, in \cite{chakraborty2022numericalradii} there is a "$-$" instead of a "$+$" in front of $2z\RE\left(A_{\alpha_{|k}}\right)$ in the definition of $Q_{k,\alpha}$. Because $S_\alpha$ is unitary equivalent to $S_{-\alpha}$, they share the same numerical range. Thus for our problem, even if it means considering $S_{-\alpha}$ instead of $S_\alpha$, we can choose to define $Q_{k,\alpha}$ this way. We do so in order to make the computations easier to follow.

The main goal of this article is to give a characterisation of the fact that $w(S_\alpha)=w_e(S_\alpha)=(a+b)/2$ which is more tractable than the previous one. Indeed, thanks to Theorem \ref{central2} and Remark \ref{recc}, we obtain an effective algorithm to apply to the $(a_i)$ to determine if $W(S_\alpha)=D\left(0,\frac{a+b}{2}\right)$ or not. The cost of this is the fact that we lose the information about the numerical radius when $w(S_\alpha)>(a+b)/2$.

The paper is organised as follows: Section \ref{TORM} is devoted to the proof of our main result Theorem \ref{central2} and Section \ref{remarques} contains applications and remarks.

\section{When $\overline{W(S_\alpha)}=\overline{D}\left(0,\frac{a+b}{2}\right)$?}\label{TORM}

We start by setting some notations. We pose $Q_0=a+b$ and $Q_{-1}=1$ and for $k>0$, we define $Q_k$ by the following determinant: 
\begin{equation*}
    Q_k:=\begin{vmatrix}
    a+b&\alpha_0&&\\
    \alpha_0&a+b&\ddots&\\
    &\ddots&\ddots&\alpha_{k-1}\\
    &&\alpha_{k-1}&a+b
\end{vmatrix}.
\end{equation*}

We define $R_0=a+b$ and $R_0^{(1)}=1$. For $1\leq k\leq n$ we set

\begin{equation*}
    R_k:=\begin{vmatrix}
    a+b&a_{n-k+1}&&\\
    a_{n-k+1}&a+b&a_{n-k}&\\
    &\ddots&\ddots&\ddots\\
    &&a_n&a+b
\end{vmatrix}\text{ and }
R_k^{(j)}:=\begin{vmatrix}
    a+b&a_{n-k+1}&&\\
    a_{n-k+1}&a+b&a_{n-k}&\\
    &\ddots&\ddots&\ddots\\
    &&a_{n-j}&a+b
\end{vmatrix} \text{ for }j\in\llbracket 0,k-1\rrbracket.
\end{equation*}
For $j=k$ and $j=m+1$, we set $R_k^{(k)}=a+b$ and $R_k^{(k+1)}=1$.

As suggested by the nomenclature, we will see in the following that the objects defined above have a link with the polynomials $Q_{\alpha,k}$ defined in the Introduction.

Finally, for a finite sequence $\alpha=(a_1,\dots,a_n)$ let us set $A_\alpha=\begin{pmatrix}
    0&&&\\
    a_1&0&&\\
    &\ddots&\ddots&\\
    &&a_n&0
\end{pmatrix}$.

The next result gives a very useful way to compute by induction the quantities above.

\begin{lem}\label{DVLP2}
    For all $k\geq1$, the following equality holds:
    \begin{equation*}
        Q_k=(a+b)Q_{k-1}-\alpha_{k-1}^2Q_{k-2}.
    \end{equation*}
    Furthermore, for all $k \in \llbracket 2,n\rrbracket$ and $j\in\llbracket 0,k-1\rrbracket$, we have
    \begin{equation*}
        R_k^{(j)}=(a+b)R_{k-1}^{(j)}-a_{n-k+1}^2R_{k-2}^{(j)}\text{ and }R_k^{(j)}=(a+b)R_k^{(j+1)}-a_{n-j}^2R_k^{(j+2)}.
    \end{equation*}
\end{lem}

\begin{proof}
    First note that $Q_1=(a+b)^2-\alpha_0^2=(a+b)Q_0-\alpha_0^2Q_{-1}$ and
    \begin{equation*}
        Q_2=(a+b)^3-(a+b)\alpha_0^2-(a+b)\alpha_1^2=(a+b)Q_1-\alpha_1^2Q_0.
    \end{equation*}
    
    Then, let us assume that $k\geq 3$. We develop the determinant $Q_k$ following the last line. We get:
    \begin{align*}
        Q_k=&\begin{vmatrix}
    a+b&\alpha_0&&\\
    \alpha_0&a+b&\ddots&\\
    &\ddots&\ddots&\alpha_{k-1}\\
    &&\alpha_{k-1}&a+b
\end{vmatrix}\\
    =&(a+b)Q_{k-1}-\alpha_{k-1}\begin{vmatrix}
    a+b&\alpha_0&&\\
    \alpha_0&a+b&\alpha_1&\\
    &\ddots&\ddots&\ddots\\
    &&\ddots&\ddots&\alpha_{k-3}&0\\
    &&&\alpha_{k-3}&a+b&0\\
    &&&&\alpha_{k-2}&\alpha_{k-1}
    \end{vmatrix}\\
    =&(a+b)Q_{k-1}-\alpha_{k-1}^2Q_{k-2}.
    \end{align*}
    The two other assertions are shown the same way.
\end{proof}

\begin{lem}\label{infini2}
    Let $a,b>0$ and $(U_k)_{k\geq0}$ a sequence such that for $k\geq2$, $U_k=(a+b)U_{k-1}-a^2U_{k-2}$ if $k$ is even and $U_{k}=(a+b)U_{k-1}-b^2U_{k-2}$ if $k$ is odd. For $k\geq0$, we have:
    \begin{equation*}
        U_{2k+1}=(ka+(k+1)b)a^kb^{k-1}U_1-k(a+b)a^{k+1}b^{k-1}U_0
    \end{equation*}
    and
    \begin{equation*}
        U_{2k+2}=(k+1)(a+b)a^kb^kU_1-((k+1)a+kb)a^{k+1}b^kU_0.
    \end{equation*}
    We also have:
    \begin{equation*}
        U_{2k+1}=a^kb^{k-1}\left[k(a+b)(U_1-aU_0)+bU_1\right]
    \end{equation*}
    and
    \begin{equation*}
        U_{2k+2}=a^kb^k\left[(k+1)(a+b)(U_1-aU_0)+abU_0\right].
    \end{equation*}
\end{lem}

\begin{proof}
    The third and fourth equalities are simple consequences of the first two ones. We will prove them by mathematical induction. For $k=0$ the first equality gives $U_1=U_1$ and the second gives $U_2=(a+b)U_1-a^2U_0$.

    Now let us assume that the two equalities hold for $k\geq 0$. By definition of the sequence $(U_k)_k$, we have:
\begin{align*}
    U_{2k+3}&=(a+b)U_{2k+2}-b^2U_{2k+1}\\
    &=\left[(k+1)(a+b)^2a^kb^k-(ka+(k+1)b)a^kb^{k+1}\right]U_1\\
    -&\left[(a+b)((k+1)a+kb)a^{k+1}b^k-k(a+b)a^{k+1}b^{k+1}\right] U_0\text{ by hypothesis}\\
    &=\left[ (k+1)a^{k+2}b^k+(k+1)a^kb^{k+2}+2(k+1)a^{k+1}b^{k+1}-ka^{k+1}b^{k+1}-(k+1)a^kb^{k+2}\right]U_1\\
    -&(a+b)\left[(k+1)a^{k+2}b^k+ka^{k+1}b^{k+1}-ka^{k+1}b^{k+1}\right]U_0\\
    &=\left[(k+1)a^{k+2}b^k+(k+2)a^{k+1}b^{k+1} \right]U_1-(k+1)(a+b)a^{k+2}b^kU_0\\
    &=((k+1)a+(k+2)b)a^{k+1}b^kU_1-(k+1)(a+b)a^{k+2}b^kU_0.
\end{align*}
On the other hand, we have:
\begin{align*}
    U_{2k+4}&=(a+b)U_{2k+3}-a^2U_{2k+2}\\
    &=(a+b)\left[((k+1)a+(k+2)b)a^{k+1}b^k-(k+1)a^{k+2}b^k\right]U_1\\
    -&\left[(a+b)^2(k+1)a^{k+2}b^k-((k+1)a+kb)a^{k+3}b^k\right]U_0\\
    &=(a+b)\left[(k+1)a^{k+2}b^k+(k+2)a^{k+1}b^{k+1}-(k+1)a^{k+2}b^k\right]U_1\\
    -&\left[(k+1)a^{k+4}b^k+(k+1)a^{k+2}b^{k+2}+2(k+1)a^{k+3}b^{k+1}-(k+1)a^{k+4}b^k-ka^{k+3}b^{k+1}\right]U_0\\
    &=(k+2)(a+b)a^{k+1}b^{k+1}U_1-\left[(k+2)a^{k+3}b^{k+1}+(k+1)a^{k+2}b^{k+2}\right]U_0\\
    &=(k+2)(a+b)a^{k+1}b^{k+1}U_1-((k+2)a+(k+1)b)a^{k+2}b^{k+1}U_0.
\end{align*}
Thus, the lemma is proved by mathematical induction.
\end{proof}

The next lemma which is a generalisation of Lemma \ref{DVLP2} will be needed in the proof of Theorem \ref{central2}.

\begin{lem}\label{les carrés}
    For all $k\in \llbracket 1,n\rrbracket$, $j \in\llbracket 0, k-1\rrbracket$ and $l \in \llbracket j, k-1\rrbracket$, we have
    \begin{equation*}
        R_k^{(j)}=R_l^{(j)}\cdot R_k^{(l+1)}-a_{n-l}^2\cdot R_{l-1}^{(j)}\cdot R_k^{(l+2)}.
    \end{equation*}
\end{lem}

\begin{proof}
    Let us fix $k\in \llbracket 1,n\rrbracket$ and $j \in\llbracket 0, k-1\rrbracket$. We will prove the lemma by mathematical induction on $l\in \llbracket j, k-1\rrbracket$.

    For $l=j$, Lemma \ref{DVLP2} gives 
    \begin{align*}
        R_k^{(j)}&=(a+b)R_k^{(j+1)}-a_{n-j}^2R_k^{(j+2)}\\
        &=R_j^{(j)}\cdot R_k^{(j+1)}-a_{n-j}^2\cdot R_{j-1}^{(j)}\cdot R_k^{(j+2)} \text{ with }R_j^{(j)}=a+b\text{ and } R_{j-1}^j=1.
    \end{align*}

    Assume now that the property is true up to an integer $l$. We then have by assumption
    \begin{align*}
        R_k&=R_l^{(j)}\cdot R_k^{(l+1)}-a_{n-l}^2\cdot R_{l-1}^{(j)}\cdot R_k^{(l+2)}\\
        &=R_l^{(j)}\left((a+b)R_k^{(l+2)}-a_{n-l-1}^2R_k^{(l+3)}\right)-a_{n-l}^2\cdot R_{l-1}^{(j)}\cdot R_k^{(l+2)} \text{ by Lemma \ref{DVLP2} applied to }R_k^{(l+1)}\\
        &=\left((a+b)R_l^{(j)}-a_{n-l}^2R_{l-1}^{(j)}\right)R_k^{(l+2)}-a_{n-l-1}^2\cdot R_l^{(j)}\cdot R_k^{(l+3)}\\
        &=R_{l+1}^{(j)}\cdot R_k^{(l+2)}-a_{n-l-1}^2\cdot R_{l}^{(j)}\cdot R_k^{(l+3)} \text{ by Lemma \ref{DVLP2}.}
    \end{align*}
    Hence, the result is proved by mathematical induction.
\end{proof}

We now give an intermediate result to reach Theorem \ref{central2} which is interesting in itself.

\begin{thm}\label{Sylvester augmenté}
    Let $\alpha\in \ell^\infty$, with our notations $w(S_\alpha)\leq (a+b)/2$ if and only if $Q_k\geq 0$ for all $k\geq1$.
\end{thm}

\begin{proof}
    Let us first prove that $w(S_\alpha)\leq (a+b)/2$ if and only if every finite minor of a certain infinite matrix is non-negative.

    We have:
    \begin{align*}
        w(S_\alpha)\leq (a+b)/2&\Leftrightarrow w(S_{-\alpha})\leq (a+b)/2 \text{ because $S_\alpha$ and $S_{-\alpha}$ are unitary equivalent}\\
        &\Leftrightarrow w(2\RE(S_{-\alpha}))\leq a+b \text{ since $W(S_{-\alpha})$ is a disk centred at 0}\\
        &\Leftrightarrow \langle \left((a+b)I-2\RE(S_{-\alpha})\right)x,x\rangle \geq0\\
        &\Leftrightarrow \text{the infinite matrix }M:=\begin{pmatrix}
            a+b&\alpha_0&&&\\
    \alpha_0&a+b&\alpha_1&&\\
    &\alpha_1&\ddots&\ddots&\\
    &&\ddots&\ddots&\ddots
        \end{pmatrix} \text{ is positive.}
    \end{align*}
    For $J\subseteq \N\cup\{0\}$ finite, we will denote $M_J$ the sub-matrix of $M$ built by the columns and lines indexed by the elements of $J$. $M\geq0$ if and only if $M_J\geq0$ for every $J\subseteq \N\cup\{0\}$ finite. Indeed $\langle M x,x\rangle\geq 0,\forall x\in \ell^2$ implies that $\langle M x,x\rangle\geq 0,\forall x$ a finite sequence, that is $M_J\geq 0$ for all $J$ finite. Reciprocally, suppose $M_J\geq 0$ for all $J$ finite and choose $x\in \ell^2$. We define $J_n=\{0,1,\dots,n\}$ and $x^{(n)}=(x_0,x_1,\dots , x_n,0,0,\dots)$. Since $M$ is bounded, we have
    \begin{equation*}
        0\leq\langle M_{J_n}x_n,x_n\rangle=\langle M x_n,x_n\rangle\to_n \langle Mx,x\rangle.
    \end{equation*}
    Thus, $\langle Mx,x\rangle\geq0$ and then $M\geq0$ and we obtain the desired equivalence.

    Sylvester's Criterion states that for a finite matrix $M_J$, $M_J\geq 0$ if and only if every minor of $M_J$ is non-negative (\textit{i.e.} $\det\left(M_K\right)\geq0$ for all $K\subseteq J$). Thus, $M\geq0$ boils down to $\det\left(M_J\right)\geq0$ for all $J\subseteq \N\cup\{0\}$ finite. To finish the proof, it remains to show that we can restrict ourselves to the $(Q_k)_{k\geq0}$, that is to the $J_k=\{0,1,\dots,k\}$.

    If $J=\{i_1,\dots, i_n\}$, $\det\left(M_J\right)=\begin{vmatrix}
        m_{i_1,i_1}&m_{i_1,i_2}&\dots&m_{i_1,i_n}\\
        m_{i_2,i_1}&m_{i_2,i_2}&\dots&m_{i_2,i_n}\\
        \vdots&\vdots&\dots&\vdots\\
        m_{i_n,i_1}&\dots&\dots&m_{i_n,i_n}
    \end{vmatrix}$, where $M=(m_{i,j})$. Since $M$ is tridiagonal, so is $M_J$ and if there are two non-consecutive integers in $J$, $M_J$ is a block-diagonal matrix. Thus, the non-negativity of each finite minor of $M$ amounts to the non-negativity of $\det(M_J)$ for all the $J$ built of consecutive integers. Let $J=\{n-l+1,n-l+2,\dots,n\}$ be a such set of integers. With the notations settled at the beginning of the section, we have $\det\left(M_j\right)=R_l$. By Lemma \ref{les carrés}, we have for all $l\in \llbracket 0, n-1\rrbracket$:
    \begin{equation}\label{sauve}
        Q_n=R_n^{(0)}=R_l^{(0)}\cdot R_n^{(l+1)}-a_{n-l}^2\cdot R_{l-1}^{(0)}\cdot R_n^{(l+2)}=R_l\cdot Q_{n-l-1}-a_{n-l}^2\cdot R_{l-1}\cdot Q_{n-l-2}.
    \end{equation}
    Assume all the $Q_k\geq0$, by Lemma \ref{DVLP2}, that implies that all the $Q_k>0$. Taking $l=1$ in the preceding equation, we obtain $Q_n=R_1\cdot Q_{n-2}-a_{n-1}^2\cdot (a+b)\cdot Q_{n-3}\geq 0$. Thus $R_1\geq0$, and Equation \ref{sauve} for $l=2$ then ensures that $R_2\geq0$. Iterating this argument leads to $R_l\geq 0$ for all $l$. Thus, $Q_k\geq 0$ for all $k\geq 0$ implies that every finite minor of $M$ is non-negative. Since the converse also holds, the proof is complete ($Q_0=a+b>0$).
\end{proof}

We can now stand our main result.

\begin{thm}\label{central2}
    Let $n$ be an even number, $a_1,\dots,a_n>0$, and let $\alpha=(a_1,\dots,a_n,a,b,a,b,\dots)$. With the same notations as above, the following are equivalent:
    \begin{enumerate}
        \item\label{interet2} $w(S_\alpha)\leq(a+b)/2$,
        \item\label{mauvais ordre2} $Q_k\geq0$ for all $k\in\llbracket 1,n\rrbracket$ and $Q_n-aQ_{n-1}\geq0$,
        \item\label{bon ordre2} $R_k-aR_k^{(1)}\geq0$ for all $k\in\llbracket 1,n\rrbracket$.
    \end{enumerate}
\end{thm}

\begin{proof}
    $\bullet~ \ref{interet2}\Leftrightarrow\ref{mauvais ordre2}$: $w(S_\alpha)\leq (a+b)/2$ if and only if $Q_k\geq0$ for all $k\geq 1$ by Theorem \ref{Sylvester augmenté}.
    
    Lemma \ref{infini2} applied to $(U_k)_k=\left(Q_{n+k-1}\right)_k$ gives
    \begin{equation*}
        Q_{n+2k}=a^kb^{k-1}\left[k(a+b)(Q_n-aQ_{n-1})+bQ_n\right]
    \end{equation*}
    and 
    \begin{equation*}
        Q_{n+2k+1}=a^kb^k\left[(k+1)(a+b)(Q_n-aQ_{n-1})+abQ_{n-1}\right],
    \end{equation*}
    for all $k\geq 0$. The desired equivalence follows since $a,b>0$.
    
    $\bullet~ \ref{bon ordre2}\Rightarrow\ref{mauvais ordre2}$: We start by showing that if \ref{bon ordre2} holds, $R_k-a  R_k^{(1)}>0$, $\forall k<n$. Let us assume that there exists $k<n$ such that $R_k-a  R_k^{(1)}=0$. Let $k_0$ be the smallest such $k$. We must have $k_0>0$ since $R_0-aR_0^{(1)}=b>0$. By assumption and Lemma \ref{DVLP2}, we have:
    \begin{equation*}
        R_{k_0+1}-a  R_{k_0+1}^{(1)}=(a+b)\left(R_{k_0}-a  R_{k_0}^{(1)}\right)-a_{n-k_0}^2\left(R_{k_0-1}-a  R_{k_0-1}^{(1)}\right)\geq 0, \text{ since }k_0<n.
    \end{equation*}
    Thus we get $R_{k_0-1}-a  R_{k_0-1}^{(1)}\leq0$ since $a_{n-k_0}^2>0$. That is in contradiction with the choice of $k_0$. We can now assume that $R_k-a  R_k^{(1)}>0$ if $k\neq n$.

    To prove the implication, first notice that $Q_n=R_n$ and $Q_{n-1}=R_n^{(1)}$, thus it suffices to show that $Q_k\geq 0$ for all $k\in\llbracket 1,n\rrbracket$. We define $B:=\begin{pmatrix}
    0&&&\\
    a_n&0&&\\
    &\ddots&\ddots&\\
    &&a_1&0
\end{pmatrix}\in \mathcal{M}_{n+1}(\C)$. $B$ is unitary equivalent to the transpose of $A_\alpha$. Thus $W(A_\alpha)$ and $W(B)$ are the same. If we set $\alpha'=(a_1,\dots,a_n,0,0,\dots)$ and $\beta'=(a_n,\dots,a_1,0,0,\dots)$, we then have $W(S_{\alpha'})=W(A_\alpha)=W(B)=W(S_{\beta'})$. Theorem \ref{Sylvester augmenté} applied to $S_{\alpha'}$ and $S_{\beta'}$ then ensures that $Q_k\geq 0$ if and only if $R_k\geq 0$ for all $k\leq n$. Let us fix $k\in\llbracket 1,n\rrbracket$. Lemma \ref{les carrés} gives that for every $l\in \llbracket 1,k-1\rrbracket$ we have
\begin{equation}\label{l'ancien lemme}
    R_k-aR_k^{(1)}=\left(R_l-aR_l^{(1)}\right)R_k^{(l+1)}-a_{n-l}^2\left(R_{l-1}-aR_{l-1}^{(1)}\right)R_k^{(l+2)}.
\end{equation}
Since $R_k-aR_k^{(1)}\geq 0$ by assumption, the preceding equation with $l=k-2$ leads to
\begin{equation*}
    \left(R_{k-2}-aR_{k-2}^{(1)}\right)R_k^{(k-1)}-a_{n-k+2}^2\left(R_{k-3}-aR_{k-3}^{(1)}\right)R_k^{(k)}\geq0.
\end{equation*}
$R_k^{(k)}=a+b>0$, so we have $R_k^{(k-1)}\geq0$. Applying this argument successively for $l=k-1,k-2,\dots ,1$ leads to $R_k^{(l)}\geq0$ for all $l\geq2$. Since $bR_k^{(1)}-a_n^2R_k^{(2)}=R_k-aR_k^{(1)}\geq0$, $R_k^{(1)}\geq 0$ and then $R_k\geq 0$.

$\bullet ~ \ref{mauvais ordre2}\Rightarrow \ref{bon ordre2}$: Suppose that $Q_k\geq 0$ for all $k\in \llbracket 1,n\rrbracket$ and that $R_n-aR_n^{(1)}=Q_n-aQ_{n-1}\geq0$. As shown in the previous point, we have that $Q_k>0$ for all $k\in \llbracket 1,n\rrbracket$. Equality \ref{l'ancien lemme} with $k=n$ leads to
\begin{equation}\label{k=n}
    0\leq R_n-aR_n^{(1)}=\left(R_l-aR_l^{(1)}\right)R_n^{(l+1)}-a_{n-l}^2\left(R_{l-1}-aR_{l-1}^{(1)}\right)R_n^{(l+2)}
\end{equation}
for all $l\in\llbracket 1,n-1\rrbracket$. Since $R_n^{(j)}=Q_{n-j}$, Equation \ref{k=n} gives
\begin{equation*}
    \left(R_l-aR_l^{(1)}\right)Q_{n-l-1}-a_{n-l}^2\left(R_{l-1}-aR_{l-1}^{(1)}\right)Q_{n-l-2}\geq0,~l\in\llbracket 1,n-1\rrbracket.
\end{equation*}
For $l=1$, we then have $\left(R_1-aR_1^{(1)}\right)Q_{n-2}-a_{n-1}^2\left(R_{0}-aR_{0}^{(1)}\right)Q_{n-3}\geq0$ and since $R_{0}-aR_{0}^{(1)}=b>0$ and the $Q_k$ are positive, we deduce that $R_1-aR_1^{(1)}>0$. Then, inducing this argument, we obtain point \ref{bon ordre2}.
\end{proof}

\begin{rem}\label{preuve originelle}
    In this remark, we will present another proof of $\ref{interet2}\Rightarrow \ref{bon ordre2}$ which uses the work done in \cite{chakraborty2022numericalradii}. In fact, that was our original argument in the proof of Theorem \ref{central2} and it gives us the intuition of our main result.
    
    In the following, we will use the notation $G_{\alpha,o}$ and $G_{\alpha,e}$ settled in the introduction.

    Let us suppose there is a $k$ such that $R_k-a  R_k^{(1)}<0$ and set $\beta=(a_{n-k+1},\dots,a_n)$. We will also denote by $\beta$ the infinite sequence $(a_{n-k+1},\dots,a_n,a,b,a,b,\dots)$ or $(a_{n-k+1},\dots,a_n,b,a,b,a,\dots)$ depending on the parity of $k$.
    
    \ding{226} If $k$ is even, the finite sequence $\beta$ has an even length and so by \cite{chakraborty2022numericalradii}, $w(S_\beta)>(a+b)/2$ if and only if $G_{\beta,e}$ has a root between 0 and 1. $G_{\beta,e}(0)=abQ_{n-k-1,\beta}(0)=ab^{(n-k+2)/2}>0$ since the sequence $\beta$ has length $n-k$. On the other hand,
    \begin{equation*}
        G_{\beta,e}(1)=abQ_{n-k-1,\beta}(1)-a_n^2\times a\times Q_{n-k-2,\beta}(1).
    \end{equation*}
    Now, $Q_{n-k-1,\beta}(1)=\det\left((a+b)I_{n-k}+2\RE\left(A_{\beta_{|n-k-1}}\right)\right)=R_k^{(1)}$ and $Q_{n-k-2,\beta}(1)=R_k^{(2)}$ with our notations. Hence,
    \begin{align*}
        G_{\beta,e}(1)&=a\left(bR_k^{(1)}-a_n^{2}R_k^{(2)}\right)\\
        &=a\left((a+b)R_k^{(1)}-a_n^2R_k^{(2)}-aR_k^{(1)}\right)\\
        &=a\left(R_k-aR_k^{(1)}\right)\text{ by Lemma \ref{DVLP2}.}
    \end{align*}
    So, by hypothesis, $G_{\beta,e}(1)<0$ and since $G_{\beta,e}(0)>0$, the real polynomial $G_{\beta,e}$ has a root between 0 and 1 and we can claim that $w(S_\beta)>(a+b)/2$.
    
    But $S_\beta$ is a compression of $S_\alpha$ and so $w(S_\alpha)\geq w(S_\beta)>(a+b)/2$.
    
    \ding{226} If $k$ is odd, the sequence $\beta$ has an odd length and $w(S_\beta)>(a+b)/2$ if and only if $\tilde{G}_{\beta,o}$ has a root between 0 and 1, where $\tilde{G}_{\beta,o}$ is the polynomial $G_{\beta,o}$ where the roles of $a$ and $b$ have been swapped. Indeed, $S_\beta$ is no longer a finite perturbation of the shift of weight sequence $(a,b,a,b,\dots)$ but a perturbation of the one with weight sequence $(b,a,b,a,\dots)$ since $k$ is odd and $n$ is even.
    
    We have $\tilde{G}_{\beta,o}(0)=ab^{(n-k+2)/2}>0$. On the other hand, we have
    \begin{align*}
        \tilde{G}_{\beta,o}(1)&=bQ_{n-k-1,\beta}(1)-a_n^2 Q_{n-k-2,\beta}(1)\\
        &=bR_k^{(1)}-a_n^2R_k^{(2)}\\
        &=R_k-aR_k^{(1)}.
    \end{align*}
    The same argument as before leads to $w(S_\alpha)\geq w(S_\beta)>(a+b)/2$.
\end{rem}

\begin{coro}\label{corol impair}
    Let $n$ be an odd number, $a_1,\dots,a_n>0$, and let $\alpha=(a_1,\dots,a_n,b,a,b,a,b,\dots)$. With the same notations as above, the following are equivalent:
    \begin{enumerate}
        \item\label{interet3} $w(S_\alpha)\leq(a+b)/2$,
        \item\label{mauvais ordre3} $Q_k\geq0$ for all $k\in\llbracket 1,n\rrbracket$ and $Q_n-bQ_{n-1}\geq0$,
        \item\label{bon ordre3} $R_k-bR_k^{(1)}\geq0$ for all $k\in\llbracket 1,n\rrbracket$.
    \end{enumerate}
\end{coro}

\begin{proof}
    We say that $\alpha=(a_1,\dots,a_n,a_{n+1},a,b,a,b,\dots)$ with $a_{n+1}=b$ and apply the preceding Theorem. We have $w(S_\alpha)\leq (a+b)/2$ if and only if $Q_k\geq 0$ for all $1\leq k\leq n+1$ and $Q_{n+1}-aQ_n\geq0$. By Lemma \ref{DVLP2}, $Q_{n+1}=b\left(Q_n-bQ_{n-1}\right)+aQ_n$. Thus, the two first points are equivalent since $b>0$.

    Above we applied Theorem \ref{central2} to $\alpha=(a_1,\dots,a_n,a_{n+1},a,b,a,b,\dots)$ with $a_{n+1}=b$. With this notation, we have $R_k-aR_k^{(1)}\geq0$ for all $k\in\llbracket 1,n+1\rrbracket$. Our goal is to show that this is the same as $R_k^{(1)}-bR_k^{(2)}$ for all $2\leq k\leq n+1$ to get the third part of the corollary's statement. This is indeed the case since $R_k-aR_k^{(1)}=b\left(R_k^{(1)}-bR_k^{(2)}\right)$ by Lemma \ref{DVLP2} and $R_1-R_1^{(1)}=ab>0$.
\end{proof}

To conclude this section, we will look at how Theorem \ref{central2} works in practice.

\begin{rem}\label{recc}
    Lemma \ref{DVLP2} gives the following useful recurrence formula: 
    \begin{equation}\label{algo}
        R_k-aR_k^{(1)}=(a+b)\left(R_{k-1}-aR_{k-1}^{(1)}\right)-a_{n-k+1}^2\left(R_{k-2}-aR_{k-2}^{(1)}\right)
    \end{equation}
    with the two first terms given by $R_1-aR_1^{(1)}=b^2+ab-a_n^2$ and $R_2-aR_2^{(1)}=b(a+b)^2-ba_{n-1}^2-(a+b)a_n^2$.
\end{rem}

\begin{ex}\label{premier ex}
    We now give an example of Theorem \ref{central2} in action on a concrete weighted shift.
t
    Let $\alpha:=\left(\sqrt{1.8},\sqrt{0.5},\sqrt{1.4},\sqrt{1.2},1,1,\dots\right)$. We want to apply Theorem \ref{central2} to determine if $w(S_\alpha)=1=w_e(S_\alpha)$. Here $n=4$ and $a=b=1$. We only have to check four inequalities:
    \begin{enumerate}
        \item $R_1-R_1^{(1)}=2-a_4^2=0.8>0$
        \item $R_2-R_2^{(1)}=4-a_3^2-2a_4^2=0.2>0$
        \item From now, we use the recurrence formula \ref{algo} given in Remark \ref{recc}.\\ $R_3-R_3^{(1)}=2\left( R_2-R_2^{(1)}\right)-a_2^2\left(R_1-R_1^{(1)}\right)=2\times0.2-0.5\times 0.8=0$.
        \item $R_4-R_4^{(1)}=2\left( R_3-R_3^{(1)}\right)-a_1^2\left(R_2-R_2^{(1)}\right)=-1.8\times0.2=-0.36<0$.
    \end{enumerate}
    Thus, we can claim that $w(S_\alpha)>1$ by Theorem \ref{central2}. 
    
    The above computations give more information than this single fact. We showed that if $\beta=\left(\sqrt{0.5},\sqrt{1.4},\sqrt{1.2},1,1,\dots\right)$, $w(S_\beta)=1$. Furthermore, since $R_3-R_3^{(1)}=0$, for all $t>0$, we have $w(S_{\beta_t})>1$ where $\beta_t:=\left(t,\sqrt{0.5},\sqrt{1.4},\sqrt{1.2},1,1,\dots\right)$.
\end{ex}

\begin{ex}
    In Theorem \ref{central2}, the inequality $R_n-aR_n^{(1)}\geq 0$ is not sufficient to have $w(S_\alpha)=1$. Indeed, if $\alpha=\left(1,1,\sqrt{1.1},\sqrt{1.5},1,1,\dots\right)$ we have $a=b=1$ and $n=4$. We use formula \ref{algo} to compute the sequence $\left(R_k-R_k^{(1)}\right)$ and we get the following:t
    \begin{center}
        \begin{tabular}{|c|c|c|}
        \hline
        Value of $k$ & Value of $a_{n-k+1}^ 2$ & Value of $R_k-R_k^{(1)}$ \\\hline
        1 &  1.5 & 0.5\\\hline
        2&1.1&-0.1\\\hline
        3&1&0.3\\\hline
        4&1&0.7\\\hline
    \end{tabular}
    \end{center}
    Since $R_2-R_2^{(1)}<0$, Theorem \ref{central2} states that $w(S_\alpha)>1$ even if $R_4-R_4^{(1)}\geq0$.
\end{ex}

\begin{rem}
    If $n=2$, by Theorem \ref{central2}, $w(S_\alpha)=1$ if and only if $R_1-aR_1^{(1)}\geq0$ and $R_2-aR_2^{(1)}\geq0$. By Lemma \ref{DVLP2}, we have
    \begin{equation}\label{les pairs suffisent}
        R_1-aR_1^{(1)}=\frac{R_2-aR_2^{(1)}+a_1^2\left(R_0-aR_0^{(1)}\right)}{a+b}.
    \end{equation}
    Thus, $w(S_\alpha)=1$ if and only if $R_2-aR_2^{(1)}\geq 0$ since $R_0-aR_0^{(1)}=b>0$ and $a_2^2>0$. As noticed in Remark \ref{recc}, $R_2-aR_2^{(1)}=b(a+b)^2-ba_{1}^2-(a+b)a_2^2$. Thus, we have the following:
    \begin{equation*}
        \text{If }n=2,~ w(S_\alpha)=1\Leftrightarrow b(a+b)^2\geq ba_{1}^2+(a+b)a_2^2.
    \end{equation*}
    This result is contained in \cite[Theorem 2.4]{chakraborty2022numerical}.
\end{rem}

Let us note that Equation \ref{les pairs suffisent} can be easily generalised for $n\geq$ and $2\leq k\leq n$:
\begin{equation*}
    R_{k-1}-aR_{k-1}^{(1)}=\frac{R_k-aR_k^{(1)}+a_{n-k+1}^2\left(R_{k-2}-aR_{k-2}^{(1)}\right)}{a+b}.
\end{equation*}
Thus, if $n$ is even and $R_k-aR_k^{(1)}\geq0$ for all $1\leq k\leq n$ even, $R_k-aR_k^{(1)}\geq0$ for all $k\leq n$ automatically. So, we have proved the following refinement of Theorem \ref{central2}:

\begin{thm}
    Let $n$ be an even number and $\alpha=\left(a_1,\dots,a_n,a,b,a,b,\dots\right)$. $w(S_\alpha)=(a+b)/2$ if and only if $R_k-aR_k^ {(1)}\geq 0$ for all $2\leq k\leq n$ even.
\end{thm}

\begin{rem}
    The preceding theorem is not a real improvement of Theorem \ref{central2}. Indeed, if we want to compute the $R_k-aR_k^{(1)}$ using the recurrence formula given in Remark \ref{recc}, we need to compute all of them, including the odd ones.
\end{rem}

Theorem \ref{central2} also enables us to compare numerical radii of two finite perturbations of 2-periodic weighted shifts, which are not perturbations of the same weighted shift.

\begin{coro}
    Let $n$ be an even number $a_1,\dots,a_n>0$, $a,b>0$ and $\alpha=(a_1,\dots,a_n,a,b,a,b,\dots)$. Assume $w(S_\alpha)=(a+b)/2$. If $c,d>0$ are such that $c+d=a+b$ and $c<a$, then $w(S_\beta)=(c+d)/2=(a+b)/2$, where $\beta=(a_1,\dots,a_n,c,d,c,d,\dots)$.
\end{coro}

\begin{ex}
    We have seen at the end of Example \ref{premier ex} that if $\alpha=\left(\sqrt{0.5},\sqrt{1.4},\sqrt{1.2},1,1,\dots\right)$, $w(S_\alpha)=w_e(S_\alpha)=1$. Thus, the above corollary ensures that $w(S_\beta)=w_e(S_\beta)=1$ if\\
    $\beta=\left(\sqrt{0.5},\sqrt{1.4},\sqrt{1.2},0.8,1.2,0.8,1.2,\dots\right)$.
\end{ex}

Keeping the same notations, one might hope that if $a+b=c+d$, $w(S_\alpha)=w(S_\beta)$ without the assumption that $w(S_\alpha)=(a+b)/2$. It is not true. For $t\in (-1,1)$, let $\gamma_t$ be the sequence $\gamma_t=\left(\sqrt{2},1-t,1+t,1-t,1+t,\dots\right)$. We know that $w(S_{\gamma_t})=1$ whenever $t\geq0$. But if $t<0$, then $b<1$ and the condition $R_2-aR_2^{(1)}\geq0$ does not hold. Indeed, $R_2-aR_2^{(1)}=b(a+b)^2-ba_{n-1}^2-(a+b)a_n^2=2b-2<0$ if $b<1$. Thus, if $t<0$, $w\left(S_{\gamma_t}\right)>1=w\left(S_{\gamma_0}\right)$ by Theorem \ref{central2}.

\section{Applications and remarks}\label{remarques}

As mentioned in the introduction, Theorem 4.3 of \cite{WangWu2011NRofWS} leads to:

\begin{thm}\label{c'est ouvert}
    Let $a_1,\dots,a_n>0$, $a,b>0$ and $\alpha=(a_1,\dots,a_n,a,b,a,b,\dots)$, if $w(S_\alpha)\leq (a+b)/2$, then $W(S_\alpha)=D(0,(a+b)/2)$.
\end{thm}

We give a direct proof of this result, more in the spirit of the proof of Theorem \ref{central2} and results in \cite{chakraborty2022numericalradii} in the appendix.

\begin{rem}
    This result does not hold with an infinite perturbation. Wang and Wu give an example of a sequence $\alpha=(\alpha_n)_n$ with $\alpha_n\to1$ such that $W(S_\alpha)=\overline{\D}$; see the remark after Theorem 3.3 in \cite{WangWu2011NRofWS}.
\end{rem}

We conclude with the study of the case where we perturb only one coordinate of a 2-periodic weighted shift. 

\begin{prop}\label{c tout seul}
    Let $\alpha$ be the sequence $(a,b,a,b,\dots)$ apart from a $c>0$ at the $m$-th position ($m\geq0)$. $W(S_\alpha)=D(0,(a+b)/2)$ if and only if $c\leq C$ with $C=a\sqrt{1+\frac{2b}{(m+2)a+mb}}$ if $m$ is even and $C=b\sqrt{1+\frac{2a}{(m+1)(a+b)}}$ if $m$ is odd.
\end{prop}

\begin{proof}
    Let us start with the case when $m$ is odd, meaning that $\alpha$ has the form $(a,b,\dots,a,c,a,b,\dots)$. We want to apply the equivalence \ref{interet2}$\Leftrightarrow$\ref{bon ordre2} of Theorem \ref{central2} to $S_\alpha$. To do so, we define the sequence $(a_1,\dots,a_{m+1})$ of length $m+1$ as $(a,b,\dots,a,c)$. We have $R_1-aR_1^{(1)}=b^2+ab-c^2$ and $R_0-aR_0^{(1)}=b$ by definition. The recurrence formula given in Remark \ref{recc} allows us to apply Lemma \ref{infini2} to $\left( R_k-aR_k^{(1)}\right)_k$. Thus, we have for all $k\leq (m-1)/2$:
    \begin{equation*}
        R_{2k+1}-aR_{2k+1}^{(1)}=a^kb^{k-1}\left[k(a+b)(b^2-c^2)+b\left(b^2+ab-c^2\right)\right]
    \end{equation*}
    and 
    \begin{equation*}
        R_{2k+2}-aR_{2k+2}^ {(1)}=a^kb^k\left[(k+1)(a+b)(b^2-c^2)+ab^2\right],
    \end{equation*}
    for all $k\leq (m-1)/2$ since $\left(R_1-aR_1^{(1)}\right)-a\left(R_0-aR_0^{(1)}\right)=b^2-c^2$. Hence, for $k\leq (m-1)/2$,
    \begin{align*}
        &R_{2k+1}-aR_{2k+1}^{(1)}\geq0\\
        \Leftrightarrow&k(a+b)(b^2-c^2)+b^3+ab^2-bc^2\geq0\\
        \Leftrightarrow&c\leq b\sqrt{\frac{(k+1)(a+b)}{ka+(k+1)b}}=b\sqrt{1+\frac{a}{ka+(k+1)b}}.
    \end{align*}
    On the other hand, for $k\leq (m-1)/2$, we have
    \begin{align*}
        &R_{2k+2}-aR_{2k+2}^ {(1)}\geq0\\
        \Leftrightarrow &(k+1)(a+b)(b^2-c^2)+ab^2\geq0\\
        \Leftrightarrow &c\leq b\sqrt{\frac{(k+2)a+(k+1)b}{(k+1)(a+b)}}=b\sqrt{1+\frac{a}{(k+1)(a+b)}}.
    \end{align*}
    Thus, since $b\sqrt{1+\frac{a}{ka+(k+1)b}}$ and $b\sqrt{1+\frac{a}{(k+1)(a+b)}}$ are decreasing when $k$ increases and\\
    $b\sqrt{1+\frac{a}{(k+1)(a+b)}}<b\sqrt{1+\frac{a}{ka+(k+1)b}}$, we get
    \begin{equation*}
        R_k-aR_k^{(1)}\geq0,\forall 1\leq k\leq m+1\Leftrightarrow c\leq b\sqrt{1+\frac{2a}{(m+1)(a+b)}}.
    \end{equation*}
    Thus, by Theorem \ref{central2}, we obtain the desired conclusion.

    Now, if $m$ is even, $\alpha$ has the form $(a,b,\dots,b,c,b,a,b,\dots)$. We want to apply Corollary \ref{corol impair}.  We have $R_0-bR_0^{(1)}=a$, $R_1-bR_1^{(1)}=(a+b)^2-c^2-b(a+b)=a^2+ab-c^2$ and so $\left(R_1-bR_1^{(1)}\right)-b\left(R_0-aR_0^{(1)}\right)=a^2-c^2$. Interchanging the roles of $a$ and $b$ in Lemma \ref{infini2} leads to
    \begin{equation*}
        R_{2k+1}-bR_{2k+1}^{(1)}=a^{k-1}b^k\left[k(a+b)(a^2-c^2)+a\left(a^2+ab-c^2\right)\right]
    \end{equation*}
    for all $k\leq m/2$ and 
    \begin{equation*}
        R_{2k+2}-bR_{2k+2}^ {(1)}=a^kb^k\left[(k+1)(a+b)(a^2-c^2)+a^2b\right],
    \end{equation*}
    for all $k\leq (m-2)/2$. Hence, for $k\leq m/2$,
    \begin{align*}
        &R_{2k+1}-bR_{2k+1}^{(1)}\geq0\\
        \Leftrightarrow & k(a+b)(a^2-c^2)+a^3+a^2b-ac^2\geq0\\
        \Leftrightarrow & c\leq a\sqrt{\frac{(k+1)(a+b)}{(k+1)a+kb}}=a\sqrt{1+\frac{b}{(k+1)a+kb}}.
    \end{align*}
    On the other hand, for every $k\leq (m-2)/2$, we have
    \begin{align*}
        &R_{2k+2}-bR_{2k+2}^ {(1)}\geq0\\
        \Leftrightarrow & (k+1)(a+b)(a^2-c^2)+a^2b\geq0\\
        \Leftrightarrow & c\leq a\sqrt{\frac{(k+1)a+(k+2)b}{(k+1)(a+b)}}=a\sqrt{1+\frac{b}{(k+1)(a+b)}}.
    \end{align*}
    Thus, $R_k-bR_k^{(1)}\geq0$ for all $0\leq k\leq m+1$ if and only if $c\leq a\sqrt{1+\frac{2b}{(m+2)a+mb}}$ and so Corollary \ref{corol impair} concludes the proof.
\end{proof}

\begin{coro}\label{c tout seul avec des 1}
    Let $\alpha$ be the sequence built of 1 everywhere apart from a $c>0$ at the $m$-th position ($m\geq0)$, then $W(S_\alpha)=\D$ if and only if $c\leq \sqrt{\frac{m+2}{m+1}}$.
\end{coro}

\begin{proof}
    If $a=b=1$, we have:
    \begin{equation*}
        a\sqrt{1+\frac{2b}{(m+2)a+mb}}=b\sqrt{1+\frac{2a}{(m+1)(a+b)}}=\sqrt{\frac{m+2}{m+1}}.
    \end{equation*}
    Then, it suffices to apply Proposition \ref{c tout seul}.
\end{proof}

The previous result is point a) of Theorem 4.6 of \cite{WangWu2011NRofWS}. The proof uses a parametrisation from Eckstein and Rácz \cite[Theorem 2.5]{eckstein1973weighted}. The principal part of the result is the following:

\begin{thm}
    Let $\alpha=(\alpha_n)_{n\in\N}\in\ell^\infty(\N,\R_+^*)$, $w(S_\alpha)\leq 1$ if and only if there exists a sequence $(t_n)_n\subseteq [-1,1)$ such that $\alpha_n^2=(1-t_n)(1+t_{n+1}),\forall n$ with $t_0=-1$.
\end{thm}

Since the proof in \cite{eckstein1973weighted} depends on an unpublished result, Wang and Wu provide a full proof of that result (\cite[Theorem 3.1]{WangWu2011NRofWS}). Using Theorem \ref{Sylvester augmenté}, we are able to give another proof of this result.

\begin{proof}
    We will use the notation $Q_n$ settled at the beginning of Section \ref{TORM} with $a=b=1$. 
    
    If $t_0,\dots,t_n$ are already built, $t_{n+1}=\frac{\alpha_n^2}{1-t_n}-1$. Since $t_n\in [-1,1)$, $t_{n+1}$ defined this way is automatically greater than -1. The existence of the sequence $(t_n)$ amounts to the fact that $t_{n+1}=\frac{\alpha_n^2}{1-t_n}-1< 1$. Let us show that $t_n=1-Q_n/Q_{n-1}$ by mathematical induction with the same notation as before.
    \begin{equation*}
        t_1=\frac{\alpha_0^2}{1-t_0}-1=1-\frac{4-\alpha_0^2}{2}=1-\frac{Q_1}{Q_0}\text{ because }t_0=-1.
    \end{equation*}
    Now assume that $t_n=1-Q_n/Q_{n-1}$. We have
    \begin{equation*}
        t_{n+1}=\frac{\alpha_n^2}{1-t_n}-1=\frac{\alpha_n^2Q_{n-1}}{Q_n}-1=1-\frac{2Q_n-\alpha_n^2Q_{n-1}}{Q_n}=1-\frac{Q_{n+1}}{Q_n} \text{ by Lemma }\ref{DVLP2}.
    \end{equation*}
    Thus $t_n< 1$ for all $n\geq 1$ is equivalent to the fact that $Q_n/Q_{n-1}>0$ for all $n$. Since $Q_0=2>0$, this is the same as asking $Q_n>0$ for all $n$. That is equivalent to $Q_n\geq0$ for all $n$ by Lemma \ref{DVLP2} since $\alpha_n>0,\forall n$. Theorem \ref{Sylvester augmenté} concludes the proof.
\end{proof}

We finish with some information about the numerical radius of the operators considered in Corollary \ref{c tout seul avec des 1}. For $c>0$ and $m\geq 0$ an integer, we will denote by $\alpha_{c,m}$ the sequence built of 1 everywhere apart from a $c>0$ at the $m$-th position. \cite[Theorem 4.6]{WangWu2011NRofWS} provides a way to compute $w\left(S_{\alpha_{c,m}}\right)$ by solving $m+1$ equations. Here, we will give some asymptotic information about this numerical radius when $t\to+\infty$ or $m\to +\infty$.

To do so, let us briefly introduce the Aluthge transform of an operator. For an operator $T$ with polar decomposition $T=U\lvert T\rvert$ where $U$ is a partial isometry with $\ker(U)=\ker(T)$, the Aluthge transform of $T$ is $\Delta(T):=\lvert T\rvert^{\frac12}U\lvert T\rvert ^\frac12$. Yamazaki's inequality for the numerical radius states that $w(T)\leq \left(\lVert T\rVert+w(\Delta(T))\right)/2$ (\cite[Theorem 2.1]{Yamazaki2007ineq}). The polar decomposition of $S_{\alpha_{c,m}}$ is $S_{\alpha_{c,m}}=SM_{\alpha_{c,m}}$, where $M_{\alpha_{c,m}}$ denotes the multiplication operator with weight $\alpha_{c,m}$ and $S$ the (unweighted) shift. Thus, $\Delta\left(S_{\alpha_{c,m}}\right)$ is the weighted shift with weight sequence\\
$(1,1,\dots,1,\sqrt{c},\sqrt{c},1,1,\dots)$, where $\sqrt{c}$ appears in $(m-1)$-th and $m$-th position if $m>0$. Otherwise, if $m=0$, $\Delta\left(S_{\alpha_{c,0}}\right)=S_{\alpha_{\sqrt{c},0}}$. Hence, for $c>1$, we get:
\begin{align*}
    w\left(S_{\alpha_{c,m}}\right)&\leq \frac{\lVert S_{\alpha_{c,m}}\rVert+w\left(\Delta\left(S_{\alpha_{c,m}}\right)\right)}{2}\\
    &\leq \frac{\lVert S_{\alpha_{c,m}}\rVert+\left\lVert\Delta\left(S_{\alpha_{c,m}}\right)\right\rVert}{2}\\
    &=\frac{c+\sqrt{c}}{2}.
\end{align*}
Thus, $w\left(S_{\alpha_{c,m}}\right)/\lVert S_{\alpha_{c,m}}\rVert\to 1/2$ when $c\to+\infty$.

We know that, $\lVert S_{\alpha_{c,m}}\rVert/2\leq w\left(S_{\alpha_{c,m}}\right)\leq \lVert S_{\alpha_{c,m}}\rVert$. The above limit can be seen as the fact that when $c$ tends to $+\infty$, $S_{\alpha_{c,m}}$ tends to have the smallest numerical range possible.

We shall now study the behaviour of $w\left(S_{\alpha_{c,m}}\right)$ when $c$ is fixed and $m\to +\infty$. Let us fix $c>1$. We have already noted in Remark \ref{preuve originelle}, that $S_{\alpha_{c,m}}$ is a compression of $S_{\alpha_{c,m+1}}$. So, the sequence $\left(w\left(S_{\alpha_{c,m}}\right)\right)_m$ is increasing and bounded by $c$. Thus, this sequence has a limit when $m\to+\infty$. Theorem 3.1 of \cite{chang2012maximizing} ensures that this limit is the numerical radius of the bilateral weighted shift with weight sequence $(\dots,1,1,c,1,1,\dots)$ where $c$ is in position 0. Theorem 4.9 of \cite{WangWu2011NRofWS} states that this numerical radius is $\left(c^2+1\right)/2c$. We summarise our conclusions in the following statement:

\begin{prop}
    Let $\alpha_{c,m}$ be the sequence built of 1 everywhere apart from a $c>0$ at the $m$-th position ($m\geq0$).
    \begin{enumerate}
        \item $W\left(S_{\alpha_{c,m}}\right)=\D$ if and only if $c\leq \sqrt{\frac{m+2}{m+1}}$. Otherwise, $W\left(S_{\alpha_{c,m}}\right)$ is a closed disk of radius strictly greater than 1.
        \item For $m\geq0$ fixed, $\lim_{c\to+\infty}w\left(S_{\alpha_{c,m}}\right)/\lVert S_{\alpha_{c,m}}\rVert= 1/2$.
        \item For $c>1$ fixed, $\lim_{m\to+\infty}w\left(S_{\alpha_{c,m}}\right)=\left(c^2+1\right)/2c$.
    \end{enumerate}
\end{prop}

\subsection*{Acknowledgments}
The author acknowledges the support of the CDP C2EMPI, together
with the French State under the France-2030 programme, the University of Lille,
the Initiative of Excellence of the University of Lille, the European Metropolis of
Lille for their funding and support of the R-CDP-24-004-C2EMPI project. The research leading to these results received funding from the French National Research Agency ANR under Grant Agreement No ANR-24-CE40-0892-02 (project COMOP).

The author also would like to thank his advisor Chafiq Benhida for his support throughout the redaction of this article. 

\subsection*{Competing interest}
None.

\subsection*{Statement about the use of AI}
No AI was used during the preparation of this work.

\bibliography{bibli.bib}
\bibliographystyle{alpha}

\appendix

\section{Openness of the numerical range when the equality between numerical radius and essential numerical radius holds}

The purpose of this appendix is to provide a proof in our context of Theorem \ref{c'est ouvert} which is different from the one of Wang and Wu. This theorem states that for $\alpha=(a_1,\dots,a_n,a,b,a,b,\dots)$, if $w(S_\alpha)=w_e(S_\alpha)$, then $W(S_\alpha)$ is open. Our proof uses Hardy space techniques and follows the same ideas as the proofs in \cite{chakraborty2022numericalradii}.

Before the proof, we recall Theorem 2 of \cite{tam1997conjecture}.

\begin{lem}\label{vp}
    Let $S_\alpha$ be a weighted shift. If $W(S_\alpha)$ is closed, $w(S_\alpha)$ is an eigenvalue of $\RE(S_\alpha)$.
\end{lem}

\begin{proof}[Proof of Theorem \ref{c'est ouvert}]
    Let $\alpha$ be as in the statement. Since $W_e(S_\alpha)=\overline{D}(0,(a+b)/2)$, if $w(S_\alpha)= (a+b)/2$, $W(S_\alpha)$ equals $D(0,(a+b)/2)$ or $\overline{D}(0,(a+b)/2)$ since $W(S_\alpha)$ is a disk. It remains to show that $W(S_\alpha)$ cannot be closed.

    If $W(S_\alpha)=\overline{D}(0,(a+b)/2)$, then $(a+b)/2$ is an eigenvalue of $\RE(S_\alpha)$ by Lemma \ref{vp}. By transposing the problem from $\ell^2$ to the Hardy space $H^2$, one can show that this is equivalent to $f$'s membership of $H^2$, where $f$ is defined by
    \begin{equation}\label{definition f}
        f(z)=\frac{G_{\alpha,e}(z)}{abz^4+(a^2+b^2-4\alpha^2)z^2+ab}
    \end{equation}
     and $\alpha=(a+b)/2$ in our case. The details are fully provided in \cite{chakraborty2022numericalradii}.
     
     With $\alpha=(a+b)/2$, the polynomial $abz^4+\left(a^2+b^2-4\alpha^2\right)z^2+ab$ has a double root at 1 and -1. Since $\left(z\mapsto \frac{1}{(1-z)^s}\right)\in H^2$ if and only if $s\in \left(0,\frac{1}{2}\right)$ (see \cite[Example 1.1.14]{martinez2007introduction}), $f$ defined as in \ref{definition f} belongs to $H^2$ only if 1 is a double root of $G_{\alpha,e}$. In this case, $f$ will belong to $H^\infty\subseteq H^2$.
     
     Notice that 1 being a root of $G_{\alpha,e}$ is the same as $R_n-aR_n^{(1)}=0$. Let us now compute $G_{\alpha,e}'(1)$. We first set
    \begin{equation*}
        f_{k,\alpha}\colon z\mapsto\sqrt{abz^4+(a^2+b^2)z^2+ab}I_{k+1}+2z\RE(A_{\alpha_{|k}}).
    \end{equation*}
    Thus, $Q_{k,\alpha}=\det\circ f_{k,\alpha}$. We denote by $\com(A)$ the cofactor matrix of $A$. So, since $D_X\det(H)=\Tr\left(^t\com(X)\cdot H\right)$, we have
    \begin{equation*}
        Q_{k,\alpha}'(z)=\Tr\left(^t\com\left(f_{k,\alpha}(z)\right)\cdot f'_{k,\alpha}(z)\right)\text{ with }f'_{k,\alpha}(z)=\frac{2abz^3+(a^2+b^2)z}{\sqrt{abz^4+(a^2+b^2)z^2+ab}}I_{k+1}+2\RE\left(A_{\alpha_{|k}}\right).
    \end{equation*}
    If we set $g_n\colon z\mapsto a_n^2\frac{(a^2z^2+ab)z^2}{\sqrt{abz^4+(a^2+b^2)z^2+ab}}$, since $G_{\alpha,e}(z)=ab Q_{n-1,\alpha}(z)-g_n(z)Q_{n-2,\alpha}(z)$, we have
    \begin{equation}\label{But equation}
        G_{\alpha,e}'(z)=abQ_{n-1,\alpha}'(z)-g_n'(z)Q_{n-2,\alpha}(z)-g_n(z)Q'_{n-2,\alpha}(z).
    \end{equation}
    We have
    \begin{align*}
        g_n'(z)&=a_n^2\frac{(4a^2z^3+2abz)\sqrt{abz^4+(a^2+b^2)z^2+ab}-(a^2z^2+ab)z^2\frac{4abz^3+2(a^2+b^2)z}{2\sqrt{abz^4+(a^2+b^2)z^2+ab}}}{abz^4+(a^2+b^2)z^2+ab}\\
        &=a_n^2\frac{(4a^2z^3+2abz)\left(abz^4+(a^2+b^2)z^2+ab\right)-(a^2z^3+abz)(2abz^4+(a^2+b^2)z^2)}{\left(abz^4+(a^2+b^2)z^2+ab\right)\sqrt{abz^4+(a^2+b^2)z^2+ab}}.
    \end{align*}
    In order to evaluate \ref{But equation} at $z=1$, we observe that $g_n(1)=a_n^2a$ and
    \begin{equation*}
        g_n'(1)=a_n^2\frac{2a(2a+b)(a+b)^2-a(a+b)^3}{(a+b)^3}=\left(a+\frac{2a^2}{a+b}\right)a_n^2.
    \end{equation*}
    On the other hand, we have $Q_{n-2,\alpha}(1)=Q_{n-2}$ with our notations and
    \begin{equation*}
        f_{k,\alpha}(1)=f_{k,\alpha}'(1)=\begin{pmatrix}
        a+b&a_1&\\
        a_1&a+b&\ddots&\\
        &\ddots&\ddots&a_k\\
        &&a_k&a+b
    \end{pmatrix}.
    \end{equation*}
    Thus, $Q'_{k,\alpha}(1)=(k+1)Q_k$ since $^t\com(A)\cdot A=\det(A)I_m$ for every $A\in\mathcal{M}_m(\C)$. Finally, we get
    \begin{equation}\label{finale}
        G'_{\alpha,e}(1)=nabQ_{n-1}-na_n^2aQ_{n-2}-\frac{2a^2}{a+b}a_n^2Q_{n-2}.
    \end{equation}
    
    Assume now that $G_{\alpha,e}(1)=G_{\alpha,e}'(1)=0$. $G_{\alpha,e}(1)=R_n-aR_n^{(1)}=0$ is the same as $bQ_{n-1}-a_n^2Q_{n-2}=0$ thanks to Lemma \ref{DVLP2} and the fact that $R_n=Q_n$ and $R_n^{(1)}=Q_{n-1}$. Hence, by \ref{finale}, $G'_{\alpha,e}(1)=0$ is the same as $Q_{n-2}=0$. Since $bQ_{n-1}-a_n^2Q_{n-2}=0$, we get $Q_{n-1}=Q_n=0$. By Lemma \ref{DVLP2}, we get $Q_k=0$ for all $k$. This is a contradiction since $Q_0=a+b>0$. Thus, $(a+b)/2$ cannot be an eigenvalue of $\RE(S_\alpha)$ and so the numerical range of $S_\alpha$ is not closed by Lemma \ref{vp}.
\end{proof}

\end{document}